\documentclass[12pt, a4paper]{article}

\usepackage{amsmath}
\usepackage{amsthm}
\usepackage{amssymb}
\usepackage{xcolor}
\usepackage{url}

\newtheorem{theorem}{Theorem}
\newtheorem{proposition}[theorem]{Proposition}

\newtheorem{remark}[theorem]{Remark}
\newtheorem{problem}[theorem]{Problem}

\newcommand{\pb}{(p)}
\newcommand{\pbm}{(p+)}

\title{Unparalleled instances of prolifickness, random walks, and square root boundaries}

\author{
Stefan Gerhold \\
TU Wien \\
\tt{sgerhold@fam.tuwien.ac.at}
\and
Friedrich Hubalek \\
TU Wien \\
\tt{fhubalek@fam.tuwien.ac.at}
}

\date{\today}

\numberwithin{equation}{section}
\numberwithin{theorem}{section}

\begin{document}

\maketitle

\begin{abstract}

We revisit the problem of influencing the sex ratio of a population by subjecting reproduction of each family to some stopping rule. As an easy consequence of the strong law of large numbers, no such modification is possible in the sense that the ratio converges to 1 almost surely, for any stopping rule that is finite almost surely.
We proceed to quantify the effects and provide limit distributions for the properly rescaled sex ratio. Besides the total ratio, which is predominantly considered in the pertinent literature, we also analyze the average sex ratio, which may converge to values different from 1.

The first part of this note is largely expository, applying  classical results and standard methods from the fluctuation theory of random walks. In the second part we apply tail asymptotics for the time at which a random walk hits a one-sided square root boundary, exhibit the differences to the corresponding two-sided problem, and use a limit law related to the empirical dispersion coefficient of a heavy-tailed distribution. 
Finally, we derive a large deviations result for a special 
stopping strategy, using  saddle point asymptotics.
\end{abstract}

% MSC2020 Vorschlaege:
% 60C05 Combinatorial probability
% 60E07 Infinitely divisible distributions; stable distributions
% 60F10 Large deviations
% 60G40 Stopping times; optimal stopping problems; gambling theory
% 60G50 Sums of independent random variables; random walks
% 41A60 Asymptotic approximations, asymptotic expansions (steepest descent, etc.) 
\smallskip

\section{Introduction}

We revisit a problem that appears already in an old article in the American Mathematical Monthly by
Robbins~\cite{Ro52}, in a classical book by Schelling \cite[p.~72]{Sch06},  as a Google job interview question, and has been the subject of many other publications and online discussions; further
references are given below.
Assume that $n$ families have children until, in each family, the first boy arrives. What is the effect
on the total sex ratio of all children?
Let $({\xi}_{j})_{j\in\mathbb N}$ be a family of i.i.d.\ Rademacher random variables, i.e.\
symmetrically distributed on $\{-1,1\}$, and
\begin{equation}\label{eq:rw}
  S_k := \sum_{j=1}^k \xi_j,\quad k\geq1,
\end{equation}
be the corresponding random walk. For a stopping time~$\tau$,
the random variable $S_\tau$
models the difference between the number of boys and girls born into a family, 
if $+1$ stands for a boy
 and $-1$ for a girl. For $k=0$, we put $S_0:=0$ as usual.
The number of
girls among the first~$k$ children is $X_k:=\tfrac12(k-S_k)$,
 the number of boys is $Y_k:=\tfrac12(k+S_k)$, and we obviously have
$S_k=Y_k-X_k$. 

Besides the ``first boy'' strategy,
other stopping times for~$S$ can be considered. Robbins~\cite{Ro52}  shows
that any integrable stopping time~$\tau$ yields $\mathbb{E}[S_\tau]=0$, by a variant of the optional sampling theorem
based on the assumption $\sup_j \mathbb{E}[\xi_j]<\infty$, which of course holds
in our setting. Thus,
if all families use such a stopping time, and we define the total sex ratio as a quotient
of expectations, then it is
\begin{equation}\label{eq:sex ratio E}
  \frac{\mathbb{E}[X_\tau]}{\mathbb{E}[Y_\tau]} = \frac{\mathbb{E}[\tau-S_{\tau}]}{\mathbb{E}[\tau+S_{\tau}]}
  =1,
\end{equation}
for any number of families.
If, on the other hand, $\tau$ is such that each family stops if it has~$p$ boys more than girls, for some fixed $p\in\mathbb N$, then obviously $\mathbb{E}[S_\tau]=p>0$.
Robbins~\cite{Ro52} observes that this strategy would lead to practical problems,  as
$\mathbb{E}[\tau]=\infty$.
In the present note, we are mainly interested in this and other stopping strategies that aim at modifying the total sex ratio.
The problem is viewed as a mathematical puzzle, to be
interpreted as a coin-tossing problem if desired, without caring about
unrealistic family sizes. These arise from heavy-tailed distributions, and
may yield numbers of children far beyond the largest ones reported in reality,
e.g.\ in \cite{gentleman}, which motivated the title of the present paper.
Besides practical fertility limits~\cite{List},
it is well known that the actual probabilities of
boy babies and girl babies are not equal, and that the sex of a baby
is not independent of its elder siblings; see e.g.~\cite{CaSt05,GrJaLa18,Se95,YaKuSr13}.
Moreover, we neglect the possibility of multiple births.
For a survey of the early literature on birth control and the sex ratio,
we refer to Goodman~\cite{Go61}.

For the ``$p$ boys more'' strategy, the ratio of girls to boys in~$n$ families, observed after all families have stopped reproducing, is a non-degenerate random variable which is strictly smaller than~$1$.
Thus, statements like ``stopping
strategies cannot affect the sex ratio'' should be read in an asymptotic sense.
We will see that the sex ratio for~$n$ families a.s.\ converges to~$1$ as $n\to\infty$, whatever
strategies are used, and
that the effect of the ``$p$ boys more'' strategy on the sex ratio is of order $1/n$.
As an even more extreme example, consider the strategy ``stop if the surplus of boys
is  at least the square root of the number of the family's children''. We will see that the effect of this strategy is asymptotically \emph{smaller} than that of ``$p$ boys more''.
For both strategies, the average family size tends to infinity. The growth is even faster for the ``square root'' strategy, and this outweighs the effect of the larger stopping
threshold on the total sex ratio.

The total sex ratio of all children is not the only quantity of interest. 
In fact, the reasons why a family would use such strategies do not concern demography,
but rather the sex ratio of this family itself.
As mentioned above, for integrable stopping times the expected number
of boys in a family equals the expected number of girls.
Still, the sex ratio $X_\tau/Y_\tau$ of a single family can have expectation smaller than~$1$.
By the strong law of large numbers (SLLN), the average sex ratio, in the sense of averaging the
ratios of all families, then converges to $\mathbb{E}[X_\tau/Y_\tau]<1$.

The underlying mathematical problems belong to the realm of fluctuation theory of random walks and are connected to questions from insurance mathematics (the authors' usual business), in particular to discrete modelling of ruin problems \cite{Di17,Ger79,Ta67}.

\section{Families and stopping rules}

The simple symmetric random walk $(S_k)_{k\geq0}$
models the difference of boys and girls  among the children of a family, that
is $S_k=Y_k-X_k$ with
\begin{align*}
  Y_k &= \#\{ 1\leq j \leq k : {\xi}_j=1\} = \tfrac12(k+S_k), \\
  X_k &= \#\{ 1\leq j \leq k : {\xi}_j=-1\} = \tfrac12(k-S_k).
\end{align*}
It is well known from the basic theory of random walks
\cite[Chapter~XII]{Fe71} that the following two strategies are a.s.\ finite:
\begin{align}
   \tau^{\pb} &:=  \inf\big\{ k\geq 1:   Y_k = p\big\}, \quad p\in\mathbb N,
  \qquad \text{``$p$ boys''}, \label{eq:m boys} \\
   \tau^{\pbm} &:= \inf\big\{ k\geq 1: S_k = p \big\}, \quad p\in\mathbb N, \qquad \text{``$p$ boys more (than girls)''}. \label{eq:m boys more} 
\end{align}
 It has been argued~\cite{Asimov} that the strategy
``stop when there are twice as many boys as girls'' should result in a sex ratio of $2:1$.
This strategy does not necessarily terminate, though:
\begin{proposition}\label{prop:tau infty}
  For
  \begin{equation}\label{eq:tau double}
    \chi     :=\inf \{ k\geq 1 : Y_k \geq  2X_k \},
  \end{equation}
  we have $\mathbb{P}[\chi=\infty]=(3-\sqrt{5})/4 \approx 0.19$.
\end{proposition}
\begin{proof}
   Note that $Y_k\geq 2X_k$ is equivalent to $S_k\geq k/3$.
  Define $\hat\xi_j:=\frac12(3\xi_{j+1}-1)$ for $j\geq1$ and 
$\hat S_k=\sum_{j=1}^k\hat\xi_j$ for $k\geq0$. Then
$(\hat\xi_j)_{j\geq1}$ is an i.i.d.\ sequence with
$\mathbb{P}[\hat\xi_j=1]=1/2$ and $\mathbb{P}[\hat\xi_j=-2]=1/2$, and
$(\hat S_k)_{k\geq0}$ is a right-continuous random walk
in the
sense of \cite[p.~21]{Spi1976} 
respectively skip-free upwards in the insurance mathematics terminology \cite[Section~2.3]{Ger79}, and it is independent of $S_1$. 
Obviously $\mathbb{P}[\tau=1]=1/2$, and for $k\geq2$ we observe
$\hat S_k=(3(S_{k+1}-S_1)-k)/2$ and
\begin{align*}
\mathbb{P}[\chi=k]&=\mathbb{P}\bigg[S_t<\frac{t}{3} \text{ for } 1\leq t < k,\ \ S_k=\frac{k}{3}\bigg]\\
&=\mathbb{P}\bigg[S_1=-1,\ \ \hat S_t<2  \text{ for }  1\leq t \leq k-2, \ \ \hat S_{k-1}=2\bigg].
\end{align*}
Next we will use the independence of $S_1$ and $\hat S_1,\ldots,\hat S_k$ and 
the classical hitting time theorem, see \cite[Theorem~2]{Ka11}. 
For $k\geq3$, and a multiple of~$3$, this yields
\begin{align*}
\mathbb{P}[\chi=k] &=\frac12\cdot\frac{2}{k-1}\mathbb{P}[\hat S_{k-1}=2] \\
&= \frac1{k-1}\mathbb{P}\bigg[S_{k-1}=\frac{k}{3}+1\bigg]=
\frac1{k-1}\binom{k-1}{\frac{k}{3}-1}2^{1-k},
\end{align*}
and zero otherwise. Thus, we conclude
\begin{align*}
\mathbb{P}[\chi=\infty]&=1-\sum_{k=1}^{\infty}\mathbb{P}[\chi=k] \\
&= 1-\left(\frac12+\sum_{j=1}^\infty
\binom{3j-1}{j-1}\frac{2^{1-3j}}{3j-1}\right)=
\frac{\sqrt5-1}{4} \approx 0.309,
\end{align*}
where we have used (15.4.12) in~\cite{DLMF} to evaluate the series.
\end{proof}
We now characterize finiteness of strategies
that stop if the surplus of boys exceeds a prescribed
function~$h$ of the number of girls, using the law of the iterated
logarithm (LIL) \cite[Section~1.9]{Gu09}. Note that ``$p$ boys more'' is of this form,
and that part~(ii) of the following proposition applies to~\eqref{eq:tau double}, with  $h(x)=x$.
\begin{proposition}\label{prop:fin}
  Let $(S_k)_{k\geq0}$ be the simple symmetric random walk, 
  and $h:\mathbb{N}_0\to [0,\infty)$.
  \begin{itemize}
    \item[(i)]
  If there is $0<c<2$ with
  $h(k)\leq c \sqrt{k \log \log k}$ for~$k$ large, then 
  \[
     \tau_h := \inf \{ k \geq1 : S_k \geq  h(X_k) \}
  \]
  is a.s.\ finite.
  \item[(ii)]
    If there is $c>2$ with $h(k)\geq  c\sqrt{ k \log \log k}$ for~$k$ large,
    then $\mathbb{P}[\tau_h=\infty]>0$.
  \end{itemize}
\end{proposition}
\begin{proof} (i)
  Define $\psi(x):= \sqrt{x \log \log x}$. By the LIL,
  \[
    \limsup_{k\to\infty} \frac{S_k}{\psi\big(\tfrac12 k + o(k)\big)}
    = \limsup_{k\to\infty} \frac{\sqrt{2}S_k}{\psi(k)} = 2 \quad \text{a.s.}
  \]
  On the event $\{S_k/k\to0\}$, which has probability~$1$ by the SLLN, we  have
  $X_k=\tfrac12(k-S_k)=\tfrac12k+o(k)$, and so
  \[
    \limsup_{k\to\infty} \frac{S_k}{h(X_k)} \geq
    \limsup_{k\to\infty} \frac{S_k}{c \psi\big(\tfrac12 k + o(k)\big)} 
    = \frac2c >1 \quad \text{a.s.}
  \]

  (ii) 
  Define
  \[
   D_{N}:=\big\{\forall k\geq N: S_k \leq c \psi(X_k)\big\}, \quad N\in\mathbb N.
   \]
   Since $X_k=\tfrac12k+o(k)$ a.s.,
    the LIL implies that $\mathbb{P}[\bigcup_{N\geq1} D_N]=0$ is impossible, and so we
   can fix~$N$ with $\mathbb{P}[D_N]>0$. Clearly, there is an integer $r< 0$ with
   $\mathbb{P}[\{S_{N}=r\} \cap D_N]>0$. 
   It now suffices to force the path of~$S$ to stay
   below zero from~$1$ to~$N$.
   By the Markov property, the event
   \[
     \Big\{ \max_{1\leq k\leq N}S_k<0,\, S_N=r \Big\} \cap D_N
   \]
   has positive probability, and~$\tau_h$ is obviously infinite on it.
\end{proof}
\begin{remark}\label{rem:finite}
Alternatively, a family could stop when $S_k\geq h(k)$ instead, i.e.\ after the surplus of boys
exceeds a function of the  number of children in the family. The proof of the proposition
shows that this amounts to replacing~$2$ by~$\sqrt{2}$ in~(i) and~(ii).
In particular, $h(x)=c\sqrt{x}$ leads to an a.s.\ finite stopping time,
for any $c>0$ (see Section~\ref{se:root}).
\end{remark}
To model an arbitrary number of families, let $({\xi}_{nj})_{n,j\in\mathbb N}$ be a family of i.i.d.\ Rademacher variables.
For every $n\in\mathbb N$, we write
\begin{equation}
  {S}_{nk} := \sum_{j=1}^k {\xi}_{nj} = Y_{nk} - X_{nk}, \quad k\geq1,
\end{equation}
for the corresponding random walk. Let ${\tau}_n$ be an a.s.\ $\mathbb{N}$-valued stopping time for~$({S}_{nk})_{k\geq1}$. 
When we say that all families use the same stopping time, we mean that
\[
  \big((\xi_{1j})_{j\geq1}, \tau_1\big)\stackrel{\mathrm d}{=}
  \big((\xi_{2j})_{j\geq1}, \tau_2\big)
  \big)\stackrel{\mathrm d}{=} \dots
\]
The number of boys resp.\ girls in the $n$-th family is
\begin{align*}
  \mathbf Y_n &:= Y_{n\tau_n} =\#\{ 1\leq k \leq {\tau}_n : {\xi}_{nk}=1\}, \\
  \mathbf X_n &:= X_{n\tau_n} = \#\{ 1\leq k \leq {\tau}_n : {\xi}_{nk}=-1\}.
\end{align*}
We will analyze  the ratio, resp.\ average ratio, of girls to boys in the first~$n$ families,
\begin{equation}\label{eq:ratio}
 R_n:=\frac{\sum\limits_{k=1}^n \mathbf X_k}{\sum\limits_{k=1}^n \mathbf Y_k}, \quad
 \bar{R}_n := \frac1n \sum_{k=1}^n \frac{\mathbf X_k}{\mathbf Y_k},
\end{equation}
as well as the fraction and average fraction
\[
  F_n := \frac{\sum_{k=1}^n \mathbf X_k}{\sum_{k=1}^n \mathbf Y_k+\sum_{k=1}^n \mathbf X_k}, \quad
  \bar{F}_n:=\frac1n\sum\limits_{k=1}^n\frac{ \mathbf X_k}{ \mathbf Y_k+ \mathbf X_k}.
\]
For all concrete strategies we consider, $R_n$ is well-defined, since $\mathbf Y_n\geq1$
for all $n\in\mathbb N$.
For arbitrary strategies, it is well-defined for sufficiently large~$n$, as
$\mathbb{P}[\mathbf Y_n=0 \ \forall n\in\mathbb N]=0$.
In statements concerning the average ratio~$\bar{R}_n$, we will assume $\mathbf X_n\leq\mathbf Y_n$. 
If all families use $\tau^{\pb}$ from~\eqref{eq:m boys}, we write $\mathbf Y_n^{\pb}$, $R_n^{\pb}$, and so on.

\section{Almost sure convergence}\label{se:as}

\begin{proposition}\label{prop:as m}
  If all families use the ``$p$ boys'' strategy~\eqref{eq:m boys},
  for fixed $p\in\mathbb N$, then
  \[
    \lim_{n\to\infty}R_n^{ \pb} = 1\quad \text{and} \quad \lim_{n\to\infty}F_n^{ \pb} = \tfrac12 \quad \text{a.s.}
  \]
  Moreover, $\bar{R}_n^{ \pb} \to 1$ 
  and $\bar{F}_n^{ \pb} \to \mathbb{E}\big[\mathbf X_1^{ \pb}/(p+\mathbf X_1^{ \pb})\big]$ a.s., 
  and  the latter limit satisfies $\lim_{p\to\infty}\mathbb{E}\big[\mathbf X_1^{ \pb}/(p+\mathbf X_1^{ \pb})\big]=\tfrac12$.
\end{proposition}
\begin{proof}
  Each~$\mathbf X_n^{\pb}$ follows a negative binomial distribution,
\[
\mathbb{P}\big[\mathbf X_n^{\pb}=j\big]=
\mathbb{P}\big[\mathbf X_1^{\pb}=j\big]=\binom{p+j-1}{j}\frac1{2^{j+p}},\quad j\geq0,
\]
with expectation $\mathbb{E}\big[\mathbf X_1^{\pb}\big]=p$. It is clear that $\mathbf Y_n^{\pb}=p$
for all $n\in\mathbb N$. By the SLLN, $\frac1n\sum_{k=1}^n \mathbf X_k^{\pb} \to p$ a.s.,
which easily yields the first three assertions.
Convergence of~$\bar{F}_n^{\pb}$ follows from the SLLN, too. Using the hypergeometric
transformation formula \cite[(15.8.1)]{DLMF}, we have
\begin{align*}
  \mathbb{E}\bigg[\frac{\mathbf X_1^{\pb}}{p+\mathbf X_1^{\pb}}\bigg]
  &= \frac{p}{2^{p+1}(p+1)}\, {}_2F_1\bigg( \genfrac{}{}{0pt}{}{p+1,p+1}{p+2} \,\bigg|\, \frac{1}{2}\bigg) \\
  &=\frac{p}{2(p+1)}\,  {}_2F_1\bigg( \genfrac{}{}{0pt}{}{1,1}{p+2} \,\bigg|\, \frac{1}{2}\bigg).
\end{align*}
This converges to~$\tfrac12$ by the asymptotic expansion \cite[(15.12.3)]{DLMF}.
\end{proof}
For $p=1$, we have
\begin{equation}\label{eq:F 1st boy} 
   \mathbb{E}\big[\mathbf X_1^{(1)}/(\mathbf Y_1^{(1)}+\mathbf X_1^{(1)})\big] =
  \mathbb{E}\big[\mathbf X_1^{(1)}/(1+\mathbf X_1^{(1)})\big]=1-\log 2\approx 0.307,
\end{equation}
which provides an example where the average fraction $\bar{F}_n$
converges to a value smaller than $\tfrac12$.
%(Strictly speaking, this follows from
%the second term in the expansion of $\mathbb{E}[\mathbf X_1/(p+\mathbf X_1)]$ for $p\to\infty$.)
Generalizing the negative binomial distribution occurring in Proposition~\ref{prop:as m}, Sheps~\cite{Sh63}
 calculates the distributions of the family size etc.\ for the more general
strategy ``stop after at least~$p_1$ boys and $p_2$ girls'', including also
an upper limit on the total family size.

To study convergence of the ratio~$R_n$ etc.\ for general stopping times, it is convenient to embed the children of all families into one random walk, as
would be done in a non-parallel computer simulation.
The idea of using a single random walk and the SLLN in the following theorem
is briefly mentioned in~\cite{Asimov} and attributed to Eugene Salamin.
\begin{theorem}\label{thm:as}
 Let $\tau_n$ be as above, i.e.\ all of them are a.s.\ $\mathbb N$-valued stopping times,
 but not necessarily equal.
 %To avoid vanishing denominators, assume $\mathbf Y_n\geq1$ a.s., $n\in\mathbb N$.
 Then
 \begin{equation}\label{eq:R F conv}
   \lim_{n\to\infty} R_n = 1 \quad \text{and} \quad
   \lim_{n\to\infty} F_n = \tfrac12 \quad \text{a.s.}
 \end{equation}
 If all families use the same stopping time, then 
 \begin{equation}\label{eq:AF AR}
   \lim_{n\to\infty} \bar{F}_n=\mathbb{E}\bigg[ \frac{\mathbf X_1}{\mathbf Y_1+\mathbf X_1} \bigg] \quad \text{and} \quad 
   \lim_{n\to\infty} \bar{R}_n = \mathbb{E}\bigg[ \frac{\mathbf X_1}{\mathbf Y_1} \bigg]\ \text{a.s.,}
 \end{equation}
 where for the last assertion we additionally assume $\mathbf X_1\leq \mathbf Y_1$ a.s.
\end{theorem}
\begin{proof}
By the renewal property of random walks at
stopping times \cite[Theorem 11.10]{Ka21}, the sequence
\[
\xi_{11},\dots,\xi_{1\tau_1},\xi_{21},\dots,\xi_{2\tau_2},\xi_{31},\dots
\]
is i.i.d.\ Rademacher.
We denote it  by $(\tilde{\xi}_j)_{j\in \mathbb N}$, 
write $\tilde{S}_k:= \sum_{j=1}^k \tilde{\xi}_j$ for its random walk,
and $\sigma_n:=\sum_{k=1}^n \tau_k$, $n\in\mathbb{N}_0$.
%the family
%\[
%  (S_{\sigma_{n-1}+1}-S_{\sigma_{n-1}}, \dots S_{\sigma_n}-S_{\sigma_{n-1}}) ,\quad n\in\mathbb N,
%\]
%agrees with
%\[
%  (S_{n,1},\dots,S_{n,\tau_n}),\quad n\in\mathbb N.
%\]
%In particular, we have
Clearly, we have
\begin{align*}
  \mathbf Y_n &= \#\{ \sigma_{n-1}< j \leq \sigma_n : \tilde{\xi}_{j}=1\}, \\
  \mathbf X_n &= \#\{  \sigma_{n-1}< j \leq \sigma_n : \tilde{\xi}_{j}=-1\}
\end{align*}
and $\tilde{S}_{\sigma_n}=\sum_{k=1}^n\mathbf Y_k-\sum_{k=1}^n\mathbf X_k$, $n\in\mathbb N$.
By the SLLN,
\[
  0=\lim_{n\to\infty} \frac{\tilde{S}_{\sigma_n}}{\sigma_n} = 
  \lim_{n\to\infty} \frac{\sum_{k=1}^n\mathbf Y_k-\sum_{k=1}^n\mathbf X_k}{\sum_{k=1}^n\mathbf Y_k+\sum_{k=1}^n\mathbf X_k} \quad \text{a.s.},
\]
i.e.\ the relative surplus of boys vanishes in the limit.
By simple properties of the strictly decreasing map $0<x\mapsto \frac{1-x}{1+x}$,
this implies that the ratio~$R_n$ from~\eqref{eq:ratio} converges to~$1$,
and the fraction of girls $F_n=R_n/(1+R_n)$ to~$\frac12$.
For equal stopping times, convergence of  $\bar{F}_n$  and $\bar{R}_n$ again
follows from the SLLN. The assumption $\mathbf X_1\leq \mathbf Y_1$ ensures that $\mathbf X_1/\mathbf Y_1$
is integrable.
\end{proof}

By the dominated convergence theorem, 
\eqref{eq:R F conv} implies that $\lim_{n\to\infty} \mathbb{E}[F_n]=\tfrac12$.
For the ``first boy'' strategy ($p=1$ in~\eqref{eq:m boys}),
this was proven in a different way in~\cite{Za},
by deducing the explicit formula
$\mathbb{E}\big[F_n^{(1)}\big] = \frac{n}{2}\big(\psi(\frac{n+2}{2})-\psi(\frac{n+1}{2})\big)$,
where~$\psi$ is the digamma function.

The hitting time theorem, which we already used in the proof
of Proposition~\ref{prop:tau infty}, implies
\begin{equation}\label{eq:p+ pmf}
  \mathbb{P}\big[\mathbf X_1^{\pbm}=j\big] =  \frac{p}{p+2j} \mathbb{P}[S_{2j+p}=p]
  =  \frac{p}{p+2j}\binom{2j+p}{j}\frac{1}{2^{2j+p}},
  \quad j\geq0.
\end{equation}
We easily conclude that for ``one boy more'', the limit of $\bar{F}_n^{(1+)}$ in~\eqref{eq:AF AR}
is $\mathbb{E}\big[\mathbf X_1^{(1+)}/(1+2\mathbf X_1^{(1+)})\big]=
1-\pi/4\approx 0.21$, which is even smaller than the limit~\eqref{eq:F 1st boy}
for the ``first boy'' strategy.
Moreover, the limit of $\bar{R}_n^{(1+)}$ in~\eqref{eq:AF AR}
is $\mathbb{E}\big[\mathbf X_1^{(1+)}/( 1+ \mathbf X_1^{(1+)})\big]=
2\log 2-1\approx 0.386$.
\begin{problem}
  Is there a stopping time for which
  $\mathbb{E}[\mathbf X_1/\mathbf Y_1]< 2\log 2-1$?
\end{problem}
The following result holds, for example, for ``$p$ boys more'' with arbitrary~$p\in\mathbb N$.
\begin{proposition}\label{prop:<1}
  Let $\tau$ be an a.s.\ $\mathbb N$-valued stopping time that favors boys
  in the sense that $S_\tau\geq0$ a.s.\ and $\mathbb{P}[S_\tau>0]>0$.
  Then the limit of~$\bar{R}_n$ in~\eqref{eq:AF AR} is smaller than~$1$,
  and the limit of~$\bar{F}_n$ is smaller than~$\tfrac12$.
\end{proposition}
\begin{proof}
  With $H:=S_\tau/\tau$, the first limit is $\mathbb{E}[(\tau-S_\tau)/(\tau+S_\tau)]=\mathbb{E}[(1-H)/(1+H)]$.
  We have $0\leq H \leq 1$, and since $(1-H)/(1+H)$ strictly decreases
  w.r.t.~$H$ on $[0,1]$, the assertion on~$\bar{R}_n$ follows.
  As for~$\bar{F}_n$, note that $\mathbb{E}[\mathbf X_1/(\mathbf Y_1+\mathbf X_1)] = (1-\mathbb{E}[H])/2$.
\end{proof}
\textcolor{blue}{
 Note that the techniques we applied in this and the preceding section appear also in insurance mathematics, in particular in discrete time and state  models of actuarial ruin theory, respectively in the fluctuation theory of skip-free upwards random walks, see for example \cite{Di17,Ger79}. 
}

\section{Limit law for the ``$p$ boys more'' strategy}\label{se:limit laws}

To assess the asymptotic effect of birth control strategies, we consider
limit distributions for $1-R_n$, where~$R_n$ is defined
in~\eqref{eq:ratio}. Since $\frac12 - F_n = \tfrac12(1-R_n)/(1+R_n)$,
and $1+R_n\to2$ a.s.\ by Theorem~\ref{thm:as}, a version of Slutsky's theorem
 \cite[Theorem~2.7~(v)]{vaVa98} then easily yields limit laws for $\tfrac12-F_n$.
 Thus, we will not consider the fraction~$F_n$ in what follows.
For the ``$p$ boys'' strategy~\eqref{eq:m boys}, the central limit theorem  implies that $\sqrt{n}\big(1-R_n^{\pb}\big)$
converges to a Gaussian distribution, and so $1-R_n^{\pb}$ is of
order $1/\sqrt{n}$, roughly. In this section, we show that it is of order $1/n$ for the ``$p$ boys more'' strategy~\eqref{eq:m boys more}.

 We have $\mathbb{P}\big[\mathbf X_1^{\pbm}=j\big] \sim \frac{p}{2\sqrt{\pi}} \frac{1}{j^{3/2}}$
by~\eqref{eq:p+ pmf} and Stirling's formula,
and so the expectation of~$\mathbf X_1^{\pbm}$
 is infinite, as already observed in~\cite{Ro52}. 
Thus, for this strategy the sex ratio
cannot be defined by $\mathbb{E}[\mathbf X_1]/\mathbb{E}[\mathbf Y_1]$, as is often done in the more applied literature (cf.~\eqref{eq:sex ratio E}).

\begin{proposition}
   For $p\in\mathbb N$ and the ``$p$ boys more'' strategy~\eqref{eq:m boys more},
   we have 
\[
  \tfrac12np(1-R_n^{\pbm})\stackrel{\mathrm d}{\to}\chi_1^2
\]
as $n\to\infty$, i.e.\ the rescaled sex ratio converges to a chi-squared distribution.
\end{proposition}
\begin{proof}
  Since the tail of~$\mathbf X_1^{\pbm}$ satisfies $\mathbb{P}\big[\mathbf X_1^{\pbm}\geq j\big]\sim p/\sqrt{\pi j}$ as 
  $j\to\infty$, see above,
  the generalized CLT  \cite[Theorem 2.2.15]{EmKlMi97}  implies
  that $\frac{1}{n^2}\sum_{k=1}^n \mathbf X_k^{\pbm}$ converges in distribution 
  to the stable distribution with density
  \begin{equation}\label{eq:stable}
  \frac{p}{2\sqrt{\pi x^3}}e^{-\frac{p^2}{4x}}, \quad x>0.
\end{equation}
  The assertion follows from the continuous mapping theorem \cite[Theorem~2.7]{Bil99}, since
  \[
  \frac{n}{p}\big(1-R_n^{\pbm}\big) = \bigg(\frac{1}{n^2}\sum_{k=1}^n \mathbf X_k^{\pbm} + \frac{p}{n}\bigg)^{-1}.
  \qedhere
\]
\end{proof}
Alternatively, the limit law can be obtained from an explicit calculation
using the probability generating function (pgf) of~$\mathbf X_1^{\pbm}$,
\[
  \mathbb{E}\Big[z^{\mathbf X_1^{\pbm}}\Big] =  \left(\frac{1-\sqrt{1-z}}{z}\right)^p =: f(z)^p,
  \quad |z|<1.
\]
Clearly,  the pgf of $\sum_{k=1}^n \mathbf X_k^{\pbm}$ is $f(z)^{pn}$.
For $\Re s\geq0$, a straightforward computation shows
\[
  \lim_{n\to\infty}\mathbb{E}\bigg[\exp\Big({-\frac{s}{n^2}}\sum_{k=1}^n \mathbf X_k^{\pbm}\Big)\bigg]
  =\lim_{n\to\infty}f(e^{-s/n^2})^{pn} = e^{-p\sqrt{s}},
\]
which is the Laplace transform of the stable distribution with 
density~\eqref{eq:stable}.
Thus, $\frac{1}{n^2}\sum_{k=1}^n \mathbf X_k^{\pbm}$ converges in distribution to that law.
This is an easy variation of a well-known result on first passage times,
see e.g.\ Exercise~5 in~\cite{La16}.

A third argument is based on Donsker's theorem. The stopping times~$\sigma_n$ from
the proof of Theorem~\ref{thm:as}
specialize to the first passage times
\[
  \sigma_n = \inf\{ k\in\mathbb{N} : \tilde{S}_k \geq pn \}
    = \inf\{ t\geq0 : \tilde{S}_t \geq pn \},   
\]
where~$(\tilde{S}_t)_{t\geq0}$ is defined as linear interpolation of~$(\tilde{S}_k)_{k\in{\mathbb N}_0}$.
 According to Donsker's theorem \cite[Theorem 2.4.20]{KaSh91},
$(n^{-1}\tilde{S}_{n^2t})_{t\geq0}$ converges weakly to a Brownian motion~$W$.
The number of girls in the $n$-th family is $\mathbf X_n^{\pbm} = \tfrac12(\sigma_n-\sigma_{n-1}-p)$,
and thus the scaled average number is
\begin{equation}\label{eq:G s}
  \frac{1}{n^2}\sum_{k=1}^n \mathbf X_k^{\pbm} = \frac{\sigma_n}{2n^2}-\frac{p}{2n}
  =\frac12 \inf\{ t\geq 0 : n^{-1}\tilde{S}_{n^2t} \geq p\}-\frac{p}{2n}.
\end{equation}
By continuity of the first passage time functional $f\mapsto  \inf \{t \geq0: f(t)\geq p\}$,
see~\cite{Wh71}, and Slutsky's theorem we conclude that $\frac{1}{n^2}\sum_{k=1}^n \mathbf X_k^{\pbm}$
converges in distribution to $\frac12 \inf\{t\geq 0: W_t =p\}$.
The latter has density~\eqref{eq:stable}; see \cite[Section 2.6.A]{KaSh91}.

\section{Limit law for the square root strategy}\label{se:root}

%Up to now we have investigated strategies such as one boy, at least~$p$ boys, or one or~$p$ boys more than girls.
In contrast to
the concrete stopping times considered so far, we now investigate a strategy that intuitively would produce \emph{substantially} more boys than girls, but unlike the strategy waiting for twice as many boys, does terminate almost surely (see Remark~\ref{rem:finite}):
\begin{equation}\label{eq:sqrt}
  \tau := \inf\big\{ k\geq 1 : S_k \geq c\sqrt{k}\big\}, \quad c>0.
\end{equation}
As mentioned in the introduction, we will show that its effect on the total
sex ratio is asymptotically \emph{smaller} than that of the ``$p$ boys more'' 
strategy~\eqref{eq:m boys more}.
If all families use~\eqref{eq:sqrt}, we have $\mathbf X_n=\tfrac12\big(\tau_n - \lceil c \sqrt{\tau_n} \rceil\big)$ girls 
and $\mathbf Y_n=\tfrac12\big(\tau_n + \lceil c \sqrt{\tau_n} \rceil\big)$ boys in the $n$-th family, and
the deviation of the sex ratio from~$1$ for the first~$n$ families is
\begin{equation}\label{eq:1-R}
  1-R_n = 
   \frac{2\sum_{k=1}^n \lceil c \sqrt{\tau_k} \rceil}{\sum_{k=1}^n \big(\tau_k+\lceil c \sqrt{\tau_k} \rceil\big)}.
\end{equation}
We will argue in Section~\ref{se:leo} that
\begin{equation}\label{eq:tail}
  \mathbb{P}[\tau >k]\sim \alpha k^{-\kappa}, \quad k\to\infty,
\end{equation}
where $\alpha>0$ and $\kappa\in(0,\frac12)$ depend on~$c$.
By~\eqref{eq:tail}, the distribution of~$\sqrt{\tau}$ is in the domain of attraction
of a $2\kappa$-stable law. Although the latter distribution has no first moment, its 
sample mean and variance are of interest in statistics,
and have been thoroughly studied. To this end, it has been proven that~\eqref{eq:tail}
implies that the sequence
\begin{equation}\label{eq:conv uv}
  (U_n,V_n):=
  \bigg(n^{-\frac1\kappa} \sum_{k=1}^n \tau_k,\ n^{-\frac{1}{2\kappa}} \sum_{k=1}^n \sqrt{\tau_k}\bigg)
  \stackrel{\mathrm d}{\to} (U,V)
\end{equation}
converges in distribution to a non-degenerate random pair $(U,V)$,
with an explicit Laplace transform.
A short and easy proof of this is found in  
\cite[Theorem~2.1]{AlLaTe10}, to which we  refer for additional
references; see also \cite[Theorem~1']{LeWoZi81}.

\begin{proposition}
  If all families use the stopping strategy~\eqref{eq:sqrt}, with the same $c>0$,
  then $n^{1/(2\kappa)}(1-R_n)$
  converges in distribution to a non-degenerate limit law.
\end{proposition}
\begin{proof}
Since $\big|\lceil c \sqrt{\tau_k} \rceil- c \sqrt{\tau_k}\big|\leq1$ and $\kappa<\tfrac12$,
\[
  Z_n :=n^{-\frac{1}{2\kappa}}  \sum_{k=1}^n \big(\lceil c \sqrt{\tau_k} \rceil
     - c \sqrt{\tau_k}\big)
\]
converges to~$0$ a.s.
By  \cite[Theorem~2.7~(v)]{vaVa98},
we conclude from~\eqref{eq:conv uv}  that
\[
(U_n,V_n,Z_n)  \stackrel{\mathrm d}{\to} (U,V,0).
\]
By~\eqref{eq:1-R} and the continuous mapping theorem,
\begin{align*}
  n^{\frac{1}{2\kappa}}(1-R_n)
  &= \frac{2cn^{-\frac{1}{2\kappa}}\sum_{k=1}^n \sqrt{\tau_k}+
  2 n^{-\frac{1}{2\kappa}}\sum_{k=1}^n \big(\lceil c \sqrt{\tau_k} \rceil- c\sqrt{\tau_k}\big)
  }{n^{-\frac{1}{\kappa}} \sum_{k=1}^n \tau_k 
     + n^{-\frac{1}{\kappa}} \sum_{k=1}^n c \sqrt{\tau_k}
       + n^{-\frac{1}{\kappa}} \sum_{k=1}^n \big(\lceil c \sqrt{\tau_k} \rceil- c\sqrt{\tau_k}\big)}
       \\ 
       &= \frac{2c V_n+2 Z_n}{U_n +  cn^{-\frac{1}{2\kappa}}V_n +  n^{-\frac{1}{2\kappa}}Z_n}
\end{align*}
converges in distribution to $2cV/U$.
\end{proof}
Thus, the three concrete strategies that we have mainly investigated
result in $1-R_n$ being roughly of the following order.
\begin{enumerate}
  \item $p$ boys~\eqref{eq:m boys}:  $1/\sqrt{n}$,
  \item $p$ boys more~\eqref{eq:m boys more}:  $1/n$,
  \item Square root strategy~\eqref{eq:sqrt}: $n^{-1/(2\kappa)}$, where $2\kappa<1$
  and~$\kappa$ depends on~$c$.
\end{enumerate}

\section{Hitting time of the square root function}\label{se:leo}

The tail asymptotic~\eqref{eq:tail} is taken from~\cite{Br67}. However,
the proofs given in that paper concern the \emph{two}-sided hitting time
\[
 \inf\big\{ k\geq 1 : |S_k| \geq c\sqrt{k}\big\}, \quad c>0,
\]
for the simple symmetric random walk, which is reduced to the corresponding
problem for Brownian motion.
This should be corrected in the definition of the hitting time in the introduction of~\cite{Br67} and before
\cite[Theorem~2]{Br67}, i.e.\ $S_n$ must be replaced by $|S_n|$. In the second line of
\cite[Section~2]{Br67}, the Brownian motion $\xi(t)$ should
also be replaced by $|\xi(t)|$.
The function~$\beta$ in \cite[Theorem~1]{Br67}, which is defined in the proof of that
theorem, is the tail exponent of
the two-sided problem, and is larger than our~$\kappa$.
\begin{theorem}
  The stopping time~\eqref{eq:sqrt} satisfies~\eqref{eq:tail}, with
  $\kappa(c) \in(0,\tfrac12)$. The function~$\kappa$ decreases and satisfies
  $\lim_{c\downarrow0}\kappa(c)=\frac12$ and $\lim_{c\to\infty}\kappa(c)=0$. 
\end{theorem}
\begin{proof}
  The proof is an adaption of the one given in~\cite{Br67}. Two non-trivial modifications are required for the one-sided hitting time, which we now describe.
  First, the two-sided estimate used to establish \cite[(3.12)]{Br67} does not
  immediately work for $a_1(t)=-\infty$. It can be replaced by applying the following
  observation to \cite[(3.11)]{Br67}:
  The cdf of the square root hitting time
  $\tau(a,b,c)$, in the notation of~\cite{Nov1971}, is an analytic function
  of~$a$, if everything else is fixed. To see this, express the cdf by Mellin inversion,
  and apply (12.7.14) und (13.8.11) in~\cite{DLMF}  to the Mellin transform \cite[(1)]{Nov1971} to obtain a bound for the integrand that justifies Mellin inversion and shows analyticity.
  
  Second, the function $\Phi$ used 
in the proof of \cite[Theorem~1]{Br67} needs to be replaced by the transform
of the one-sided  exit time. By \cite[Theorem~1]{Nov1971}, this is
$e^{-c^2/4} D_{-\lambda}(0)/D_{-\lambda}(-c)$,
where~$D$ is the parabolic cylinder function \cite[\S 12.1]{DLMF}.
This transform has a zero at $\lambda=-1$, and a largest negative
pole $\lambda_0(c)\in(-1,0)$ which depends on~$c>0$.
We define $\kappa(c)=-\lambda_0(c)/2$, and the asymptotic statement~\eqref{eq:tail}
follows as in~\cite{Br67}, to which we add two more remarks.
  First, note that the case that there exists $k\in\mathbb N$ with $\mathbb{P}[\tau>k]=0$, mentioned in \cite[Theorem~2]{Br67}, cannot occur for the simple symmetric random walk.
  Second, the last line of \cite[(3.15)]{Br67}, which would imply the incorrect statement
 $\lim_{u\to\infty} Q_{n,m}(\gamma,\eta)=0$, contains an error:
The term $c_N \rho_N$ needs to be replaced
by a term  of order $\sqrt{\rho_N}$. This is incidental,
as the proof of \cite[Proposition~1]{Br67} is already complete after \cite[(3.14)]{Br67}.

  As for the  properties of~$\kappa$ we  claim, note that the function decreases because otherwise~\eqref{eq:tail}
  could not hold. The limiting values follow from \cite[(12.4.1)]{DLMF}
  and \cite[(11.3.24)]{Te15}.
\end{proof}

For the reader's convenience we list a few further minor typos in~\cite{Br67}: $\lim_{c\to\infty}$ in part~(ii) of Theorem~1 should be
$\lim_{c\to0}$, and $c^{2n}$ in the second two
lines of~(2.7) should be $2^nc^{2n}$. After~(3.13), $\lambda$ should be defined to be
$1/(e^{2u_0}-1)$, because this is the reciprocal of the true limit in~(3.13).
In~(3.14), the first $\xi(t)$ should be $\xi(1)$. For
further references on Brownian motion hitting a square root boundary, we refer to~\cite{DHSW2022}.

\section{Large deviations}

We have shown in Theorem~\ref{thm:as} that the sex ratio~$R_n$ converges
to~$1$ a.s., regardless of which stopping times are used. Large deviations concern the probability that the ratio stays away
from the limit, i.e.\ $R_n\leq 1-\varepsilon$, for fixed $\varepsilon>0$.
%One way to view large deviations statements is as a refinement of the SLLN:
%``Atypical'' events are sufficiently unlikely to infer a.s.\ convergence
%by the Borel-Cantelli Lemma (see \cite[Section~I.2]{Ho20}).
For the ``$p$ boys'' strategy~\eqref{eq:m boys}, it follows from Cram\'er's classical
theorem \cite[Theorem 2.2.3]{DeZe98} that the ratio~$R_n^{\pb}$ will deviate from~$1$ only with exponentially small
probability. We will show exponential decay for the ``$p$ boys more'' strategy~\eqref{eq:m boys more}. Here, large
deviations of the ratio amount to estimating
$\mathbb{P}[\frac1n \sum_{k=1}^n{\mathbf X_k^{\pbm}} \leq c ]$, $c>0$, for large~$n$, since
\[
  \mathbb{P}[R_n^{\pbm} \leq 1-\varepsilon]
  = \mathbb{P}\bigg[ \frac1n \sum_{k=1}^n \mathbf X_k^{\pbm} \leq p \Big(\frac{1}{\varepsilon}-1\Big)\bigg],
  \quad \varepsilon \in (0,1).
\]
Thus, we consider the probability that the
average family size stays bounded.
Note that the standard reference for large deviations of
sums of heavy tailed random variables, Vinogradov~\cite{Vi94},
deals with upper large deviations, whereas we
are interested in lower large deviations.
In the light of the end of Section~\ref{se:limit laws}, a natural way to proceed 
would be to apply the large deviations result accompanying Donsker's theorem,
which is Mogulskii's theorem \cite[Theorem~5.1.2]{DeZe98}.
However, Mogulskii's theorem concerns sample path large deviations, far more
general than what we need, and is stated for a finite time interval. As we are not
aware of an extension to the half-line, we instead give a direct estimate, by applying
the saddle point method to the probability generating function of~$\mathbf X_1^{\pbm}$.

\begin{theorem}
  For $c>0$ and $p\in\mathbb{N}$, define
  \[
    \rho(p,c) := \bigg( \frac{p+2c}{2(p+c)} \bigg)^p
    \bigg( \frac{(p+2c)^2}{4c(p+c)} \bigg)^c.
  \]
  Then $\rho(p,c)\in(0,1)$, and the average number of girls for the ``$p$ boys
  more'' strategy~\eqref{eq:m boys more} satisfies the large deviations estimate
  \[
   \mathbb{P}\bigg[\frac1n \sum_{k=1}^n{\mathbf X_k^{\pbm}} \leq c \bigg] = \rho(p,c)^{n+o(\log n)},
   \quad \ n\to\infty.
  \]
\end{theorem}
\begin{proof}
   Using the pgf $f(z)^{pn}$ from Section~\ref{se:limit laws} and Cauchy's integral formula, we calculate
   \begin{align*}
      \mathbb{P}\bigg[ \sum_{k=1}^n{\mathbf X_k^{\pbm}} \leq c n \bigg]
     &= 
      \sum_{k=0}^{\lfloor cn \rfloor}  \frac{1}{2\pi i} \oint
     \frac{f(z)^{pn}}{z^{k+1}} dz \\
     &=\frac{1}{2\pi i} \oint f(z)^{pn} \frac{z^{-\lfloor cn \rfloor-1}-1}{1-z}dz.
   \end{align*}
   First, we replace $\lfloor cn \rfloor$ by $cn$. It is easy to see that this misses
   only an oscillating factor bounded between two positive constants, which
   is not relevant at the desired asymptotic accuracy. Similarly, we can neglect
   the factor $1/(1-z)$, as long as we consider a fixed integration circle of radius
   $<1$, and $-1$ in the numerator can be removed by a straightforward bound.
   We are thus left with the integral
   \begin{equation}\label{eq:sp int}
     \frac{1}{2\pi i} \int_{|z|=\hat{z}} e^{-n \eta(z)}dz,
   \end{equation}
   where $\eta(z):=-p \log f(z)+c \log z$, and we have already moved the integration
   circle to the saddle point
   $\hat{z}:=4c(p+c)/(p+2c)^2$
   of the integrand, which satisfies $\eta'(\hat z)=0$. The saddle point method
   now proceeds by integrating the local expansion
   $\eta(z) = \eta(\hat z)+ \frac12 \eta''(\hat z)(z-\hat{z})^2+\dots$
   over a suitable part of the contour close to the saddle point. Because of the simple
   form of our integrand, it suffices to verify the conditions of Theorem~4.7.1
   in~\cite{Ol74}. The only non-obvious one is condition~(v), which asserts
   that the absolute value of the integrand must be strictly larger
   at the saddle point than on the rest of the integration contour.
   Clearly, it suffices to show this for the function $1-\sqrt{1-z}$.
   With 
   \[
     g_r(u) := 1+u-\sqrt{1-r^2+2u+u^2},\quad u\geq0,
  \]  
  and $z=re^{i p}$, $r\in(0,1)$, $p\in(-\pi,\pi]$, an elementary calculation yields
  \[
    \big|1-\sqrt{1-z}\big|^2=g_r(\sqrt{1-2r\cos p+r^2}).
  \]
  We thus need to show that~$g_r$ decreases from $1-r$ to $1+r$, 
  corresponding to $p=0$ and $p=\pi$, with a strict
  maximum at $1-r$. This is true,
  as it is easy to see that $g_r'$ has no zero, and $g_r'(1-r)<0$.
  Then, Theorem~4.7.1 in~\cite{Ol74} implies that~\eqref{eq:sp int}
  is $\exp(-n \eta(\hat z)+o(\log n))$, and actually 
  gives lower order factors and a full asymptotic expansion, if desired.
  Finally, $\rho<1$ must hold since we are approximating a probability.
\end{proof}

\section*{Acknowledgement}

We thank  Stephan Pirringer, Uwe Schmock, Ulrike Schneider and Sandra Trenovatz for helpful comments and suggestions.

\bibliographystyle{siam}
\bibliography{literature}

\begin{thebibliography}{10}

\bibitem{List}
{\em List of people with the most children}.
\newblock \url{https://en.wikipedia.org/wiki/List_of_people_with_the_most_children} [Accessed: 2025-03-14].

\bibitem{AlLaTe10}
{\sc H.~Albrecher, S.~A. Ladoucette, and J.~L. Teugels}, {\em Asymptotics of the sample coefficient of variation and the sample dispersion}, J. Statist. Plann. Inference, 140 (2010), pp.~358--368.

\bibitem{Asimov}
{\sc D.~Asimov}, {\em Answer to ``{G}oogle question: In a country in which people only want boys''}.
\newblock MathOverflow.
\newblock \url{https://mathoverflow.net/q/31066}.

\bibitem{Bil99}
{\sc P.~Billingsley}, {\em Convergence of probability measures}, John Wiley \& Sons, Inc., New York, second~ed., 1999.

\bibitem{Br67}
{\sc L.~Breiman}, {\em First exit times from a square root boundary}, in Proc. {F}ifth {B}erkeley {S}ympos. {M}ath. {S}tatist. and {P}robability ({B}erkeley, {C}alif., 1965/66), {V}ol. {II}: {C}ontributions to {P}robability {T}heory, Univ. California Press, Berkeley, CA, 1967, pp.~9--16.

\bibitem{CaSt05}
{\sc M.~A. Carlton and W.~D. Stansfield}, {\em Making babies by the flip of a coin?}, Amer. Statist., 59 (2005), pp.~180--182.

\bibitem{DeZe98}
{\sc A.~Dembo and O.~Zeitouni}, {\em Large deviations techniques and applications}, vol.~38 of Applications of Mathematics (New York), Springer-Verlag, New York, second~ed., 1998.

\bibitem{DHSW2022}
{\sc D.~E. Denisov, G.~Hinrichs, A.~I. Sakhanenko, and V.~I. Wachtel}, {\em Crossing an asymptotically square-root boundary by the {B}rownian motion}, Proc. Steklov Inst. Math., 316 (2022), pp.~105--120.

\bibitem{Di17}
{\sc D.~C.~M. Dickson}, {\em Insurance risk and ruin}, International Series on Actuarial Science, Cambridge University Press, Cambridge, second~ed., 2017.

\bibitem{DLMF}
{\em {\it NIST Digital Library of Mathematical Functions}}.
\newblock \url{https://dlmf.nist.gov/}, Release 1.2.3 of 2024-12-15.
\newblock F.~W.~J. Olver, A.~B. {Olde Daalhuis}, D.~W. Lozier, B.~I. Schneider, R.~F. Boisvert, C.~W. Clark, B.~R. Miller, B.~V. Saunders, H.~S. Cohl, and M.~A. McClain, eds.

\bibitem{EmKlMi97}
{\sc P.~Embrechts, C.~Kl\"uppelberg, and T.~Mikosch}, {\em Modelling extremal events}, Springer-Verlag, Berlin, 1997.

\bibitem{Fe71}
{\sc W.~Feller}, {\em An introduction to probability theory and its applications. {V}ol. {II}}, John Wiley \& Sons, Inc., New York-London-Sydney, second~ed., 1971.

\bibitem{Ger79}
{\sc H.~U. Gerber}, {\em An introduction to mathematical risk theory}, University of Pennsylvania, Wharton School, S.S. Huebner Foundation for Insurance Education, Philadelphia, PA, 1979.

\bibitem{Go61}
{\sc L.~A. Goodman}, {\em Some possible effects of birth control on the human sex ratio}, Ann. Human Genetics, 25 (1961), pp.~75--81.

\bibitem{GrJaLa18}
{\sc V.~Grech, W.~H. James, and J.~Lauri}, {\em On stopping rules and the sex ratio at birth}, Early Human Development, 127 (2018), pp.~15--20.

\bibitem{Gu09}
{\sc A.~Gut}, {\em Stopped random walks}, Springer Series in Operations Research and Financial Engineering, Springer, New York, second~ed., 2009.

\bibitem{Ka11}
{\sc W.~Kager}, {\em The hitting time theorem revisited}, Amer. Math. Monthly, 118 (2011), pp.~735--737.

\bibitem{Ka21}
{\sc O.~Kallenberg}, {\em Foundations of modern probability}, vol.~99 of Probability Theory and Stochastic Modelling, Springer, Cham, third~ed., 2021.

\bibitem{KaSh91}
{\sc I.~Karatzas and S.~E. Shreve}, {\em Brownian motion and stochastic calculus}, Springer-Verlag, New York, second~ed., 1991.

\bibitem{La16}
{\sc S.~Lalley}, {\em Lecture notes on one-dimensional random walks}, 2016.
\newblock \url{https://galton.uchicago.edu/~lalley/Courses/312/}.

\bibitem{LeWoZi81}
{\sc R.~LePage, M.~Woodroofe, and J.~Zinn}, {\em Convergence to a stable distribution via order statistics}, Ann. Probab., 9 (1981), pp.~624--632.

\bibitem{Nov1971}
{\sc A.~A. Novikov}, {\em On stopping times for a {Wiener} process}, Theory Probab. Appl., 16 (1971), pp.~449--456.

\bibitem{Ol74}
{\sc F.~W.~J. Olver}, {\em Asymptotics and special functions}, Computer Science and Applied Mathematics, Academic Press, New York-London, 1974.

\bibitem{Ro52}
{\sc H.~Robbins}, {\em A note on gambling systems and birth statistics}, Amer. Math. Monthly, 59 (1952), pp.~685--686.

\bibitem{Sch06}
{\sc T.~C. Schelling}, {\em Micromotives and macrobehavior: with a new preface and the Nobel lecture}, Norton, New York, 2006.

\bibitem{Se95}
{\sc C.~Seidl}, {\em The desire for a son is the father of many daughters}, Journal of Population Economics, 8 (1995), p.~185–203.

\bibitem{Sh63}
{\sc M.~C. Sheps}, {\em Effects on family size and sex ratio of preferences regarding the sex of children}, Population Studies, 17 (1963), pp.~66--72.

\bibitem{Spi1976}
{\sc F.~Spitzer}, {\em Principles of random walk}, Springer-Verlag, New York-Heidelberg, second~ed., 1976.

\bibitem{Ta67}
{\sc L.~Tak\'acs}, {\em Combinatorial methods in the theory of stochastic processes}, John Wiley \& Sons, Inc., New York-London-Sydney, 1967.

\bibitem{Te15}
{\sc N.~M. Temme}, {\em Asymptotic methods for integrals}, vol.~6 of Series in Analysis, World Scientific Publishing Co. Pte. Ltd., Hackensack, NJ, 2015.

\bibitem{vaVa98}
{\sc A.~W. van~der Vaart}, {\em Asymptotic statistics}, vol.~3 of Cambridge Series in Statistical and Probabilistic Mathematics, Cambridge University Press, Cambridge, 1998.

\bibitem{Vi94}
{\sc V.~Vinogradov}, {\em Refined large deviation limit theorems}, vol.~315 of Pitman Research Notes in Mathematics Series, Longman Scientific \& Technical, Harlow; copublished in the United States with John Wiley \& Sons, Inc., New York, 1994.

\bibitem{Wh71}
{\sc W.~Whitt}, {\em Weak convergence of first passage time processes}, J. Appl. Probability, 8 (1971), pp.~417--422.

\bibitem{gentleman}
{\sc X.~Y.}, {\em Unparalleled instance of prolifickness}, The Gentleman’s Magazine, 52, part 2 (1783), p.~753.

\bibitem{YaKuSr13}
{\sc R.~C. Yadava, A.~Kumar, and U.~Srivastava}, {\em Sex ratio at birth: a model based approach}, Math. Social Sci., 65 (2013), pp.~36--39.

\bibitem{Za}
{\sc D.~Zare}, {\em Answer to ``{G}oogle question: In a country in which people only want boys''}.
\newblock MathOverflow.
\newblock \url{https://mathoverflow.net/q/17963}.

\end{thebibliography}

\end{document}